\documentclass{book}
\newdimen\footheight
\makeatletter
\newcommand\vpt   {\edef\f@size{\@vpt}\rm}
\newcommand\vipt  {\edef\f@size{\@vipt}\rm}
\newcommand\viipt {\edef\f@size{\@viipt}\rm}
\newcommand\viiipt{\edef\f@size{\@viiipt}\rm}
\newcommand\ixpt  {\edef\f@size{\@ixpt}\rm}
\newcommand\xpt   {\edef\f@size{\@xpt}\rm}
\newcommand\xipt  {\edef\f@size{\@xipt}\rm}
\newcommand\xiipt {\edef\f@size{\@xiipt}\rm}
\newcommand\xivpt {\edef\f@size{\@xivpt}\rm}
\newcommand\xviipt{\edef\f@size{\@xviipt}\rm}
\newcommand\xxpt  {\edef\f@size{\@xxpt}\rm}
\newcommand\xxvpt {\edef\f@size{\@xxvpt}\rm}
\makeatother
\usepackage{rmh,makeidx,psfig}
\usepackage{amsmath,amsfonts,amssymb,graphicx,subfigure,url}
\newfont{\eurbx}{eurb10 at 10pt}
\newfont{\eurbvii}{eurb7 at 7pt}
\newfont{\eurbv}{eurb5 at 5pt}
\newfont{\cmmibx}{cmmib10 at 10pt}
\newfont{\cmmibvii}{cmmib10 at 7pt}
\newfont{\cmmibv}{cmmib10 at 5pt}
\newfont{\eurmx}{eurm10 at 10pt}
\newfont{\eurmvii}{eurm7 at 7pt}
\newfont{\eurmv}{eurm5 at 5pt}
\newfont{\msbmx}{msbm10 at 10pt}
\newfont{\msbmvii}{msbm10 at 7pt}
\newfont{\msbmv}{msbm10 at 5pt}
\newfont{\cmrx}{cmr10 at 10pt}
\newfont{\cmrvii}{cmr7 at 7pt}
\newfont{\cmrv}{cmr5 at 5pt}
\newfont{\cmbxx}{cmbx10 at 10pt}
\newfont{\cmbxvii}{cmbx7 at 7pt}
\newfont{\cmbxiv}{cmbx5 at 5pt}

\let\realA=\A
\makeatletter
\renewcommand\A[1]{\realA{#1}\edef\@currentlabel{\thechapter.\theacount}}
\makeatother
\def\dim{\hbox{\cmrx dim}}

\def\R{{\mathbb R}}
\def\Z{{\mathbb Z}}

\def\Q{{\mathbb Q}}
\def\vol{{\mathrm{vol \,} } }

\begin{document}
\makeatletter

\newcommand{\conj}[2]{\stepcounter{thm}
\par\addvspace{12pt}
\noindent{%
\edef\@currentlabel{\thechapter.\theacount.\thethm}
{\hvbxi CONJECTURE \thechapter.\theacount.\thethm}\hspace*{1em}{\trmitxi
#1}\par\addvspace{3pt}\noindent{\trmitx #2}\par\addvspace{3pt}}}

\newcommand{\lemma}[2]{\stepcounter{thm}
\par\addvspace{12pt}
\noindent{%
\edef\@currentlabel{\thechapter.\theacount.\thethm}
{\hvbxi LEMMA \thechapter.\theacount.\thethm}\hspace*{1em}{\trmitxi
#1}\par\addvspace{3pt}\noindent{\trmitx #2}\par\addvspace{3pt}}}

\newcommand{\theorem}[2]{\stepcounter{thm}
\par\addvspace{12pt}
\noindent{%
\edef\@currentlabel{\thechapter.\theacount.\thethm}
{\hvbxi THEOREM \thechapter.\theacount.\thethm}\hspace*{1em}{\trmitxi
#1}\par\addvspace{3pt}\noindent{\trmitx #2}\par\addvspace{3pt}}}

\newcommand{\corollary}[2]{\stepcounter{thm}
\par\addvspace{12pt}
\noindent{%
\edef\@currentlabel{\thechapter.\theacount.\thethm}
{\hvbxi COROLLARY \thechapter.\theacount.\thethm}\hspace*{1em}{\trmitxi
#1}\par\addvspace{3pt}\noindent{\trmitx #2}\par\addvspace{3pt}}}

%%%%%%%%%%%%%%%%%%%%%%%%%%
\setcounter{chapter}{23}

\runhds{M. Kahle}
{Chapter 23: Random Simplicial Complexes}

\vspace{-1pc}

 \noindent\hskip-3.8pc{\hvbxxiv 23\hskip27pt RANDOM SIMPLICIAL COMPLEXES}

\vspace{1pc}

\noindent{\hvxiv Matthew Kahle}

\vspace{4pc}

\textheight=50pc

\newcommand{\G}{\mathcal{G}}
\newcommand{\link}{\mbox{lk}}
\newcommand{\Pois}{\mbox{Pois}}
\newcommand{\prob}{\mathbb{P}}
\newcommand{\expect}{\mathbb{E}}

\Bnn{INTRODUCTION}

\noindent
Random shapes arise naturally in many contexts. The topological and geometric structure of such objects is interesting for its own sake, and also for applications. In physics, for example, such objects arise naturally in quantum gravity, in material science, and in other settings. Stochastic topology may also be considered as a null hypothesis for topological data analysis. %One might also think of random simplicial complexes as higher-dimensional analogues of random graphs. For example, certain kinds of random simplicial complexes have higher-dimensional expander properties.

In this chapter we overview combinatorial aspects of stochastic topology. We focus on the topological and geometric properties of random simplicial complexes. We introduce a few of the fundamental models in Section 23.1. We review high-dimensional expander-like properties of random complexes in Section 23.2. We discuss threshold behavior and phase transitions in Section 23.3, and Betti numbers and persistent homology in Section 23.4.

\A{MODELS}
\noindent
We briefly introduce a few of the most commonly studied models.

\B{ERD\H{O}S--R\'ENYI-INSPIRED MODELS}
\noindent
A few of the models that have been studied are high-dimensional analogues of the Erd\H{o}s--R\'enyi random graph.

\C{The Erd\H{o}s--R\'enyi random graph}
\noindent
The Erd\H{o}s--R\'enyi random graph $G(n,p)$ is the probability distribution on all graphs on vertex set $[n]=\{ 1, 2, \dots, n \}$, where every edge is included with probability $p$ jointly independently. Standard references include \cite{Bollobas} and \cite{JLR00}.

One often thinks of $p$ as a function of $n$ and studies the asymptotic properties of $G(n,p)$ as $n \to \infty$. We say that an event happens \emph{with high probability (w.h.p.)} if the probability approaches $1$ as $n \to \infty$.

Erd\H{o}s--R\'enyi showed that $\bar{p} = \log n / n$ is a sharp threshold for connectivity. In other words,: for every fixed $\epsilon > 0$, if $p \ge (1 + \epsilon) \bar{p}$ then w.h.p.\ $G(n,p)$ is connected, and if $p \le ( 1 - \epsilon) \bar{p}$ then w.h.p.\ it is disconnected. A slightly sharper statement is given in the following section. Several thresholds for topological properties of $G(n,p)$ are summarized in Table \ref{tab:gnp}.

\clearpage

%It is often fruitful to think of $G(n,p)$ as a continuous-time stochastic process, where the edges are added one at a time. For every edge $e$ in the complete graph $K_n$ ,we assign a uniform random number $x_e \in [0,1]$, independently. Then at time $t$, we the graph has vertex set $[n]$ and all edges $\{ e \mid x_e \le t \}$. At time $0$ there are no edges, and at time $1$ we have a complete graph.

\begin{table}%[htb]
\viiipt
\baselineskip=11pt
\renewcommand{\arraystretch}{.9}
\Table{Topological thresholds for $G=G(n,p)$. The column SHARP indicates whether the threshold is sharp, coarse, or one-sided sharp. The column TIGHT indicates whether there is any room for improvement on the present bound. \label{tab:gnp}}
{\begin{tabular}{| p{4cm}| l | l | c | c |}
    \hline
\rule[-4pt]{0pt}{13pt}{\hvbviii PROPERTY}
    & {\hvbviii THRESHOLD}
    & {\hvbviii SHARP}
    & {\hvbviii TIGHT}
    & {\hvbviii SOURCE}
    \\ \hline
%\rule[0pt]{0pt}{9pt}
$G$ is not $0$-collapsible
    & $ 1 / n$
    & one-sided
    & yes
    & \cite{Pittel88}
    \\
$H_1(G) \neq 0$
	&$ 1 / n$
    & one-sided
    & yes
    & \cite{Pittel88}
    \\
$G$ is not planar
	&$ 1 / n$
    & sharp
    & yes
    & \cite{LPW94}
    \\
$G$ contains arbitrary minors
	&$ 1 / n$
    & sharp
    & yes
    & \cite{AKS79}
    \\
$G$ is pure $1$-dimensional
    & $ \log n  /n$
    & sharp
    & yes
    & \cite{ER59}
    \\
$G$ is connected
    & $ \log n  /n$
    & sharp
    & yes
    & \cite{ER59}   
    \\
%$G$ is an expander
%    & $ \log n  /n$
%    & sharp
%    & yes
%    & \cite{HKP12}
%    \\
   \hline
\end{tabular}}
\renewcommand{\arraystretch}{1}
\xpt
\baselineskip=12pt
\end{table}

\vspace{-1pc}

\Bnn{GLOSSARY}

\begin{gllist}

\item {\index{Threshold}\trmbitx Threshold function:} Let $\mathcal{P}$ be a graph property. We say that $f$ is a threshold function for property $\mathcal{P}$ in the random graph $G = G(n,p)$ if whenever $p = \omega(f)$, $G$ has property $\mathcal{P}$ w.h.p.\ and whenever $p = o(f)$, $G$ does not have property $\mathcal{P}$.

\item {\index{Sharp}\trmbitx Sharp threshold:} We say that $f$ is a sharp threshold for graph property $\mathcal{P}$ if there exists a function $g = o(f)$ such that for $p < f - g$, $G \notin \mathcal{P}$ w.h.p. and if $p > f + g$, $G \in \mathcal{P}$ w.h.p.

\item {\index{simplicial complex}\trmbitx Simplicial complex:}\quad A simplicial complex $\Delta$ is a collection of subsets of a set $S$, such that (1) if $U \subset V$ is nonempty and $V \in \Delta$ then $U \in \Delta$, and (2) $\{ v \} \in \Delta$ for every $v \in S$. An element of $\Delta$ is called a face. Such a set system can be naturally associated a topological space by considering every set of size $k$ in $\Delta$ to represent a $k-1$-dimensional simplex, homeomorphic to a closed Euclidean ball. This topological space is sometimes called the geometric realization of $\Delta$, but we will slightly abuse notation and identify a simplicial complex with its geometric realization.
%The underlying set $S$ is sometimes called the ground set. The elements of $S$ represent $0$-dimensional faces, or vertices of $\Delta$.

\item {\index{link}\trmbitx Link:}\quad Given a simplicial complex $\Delta$ and a face $\sigma \in \Delta$, the link of $\sigma$ in $\Delta$ is defined by
$$ \link_{\Delta} ( \sigma ) = \{ \tau \in \Delta \mid \tau \cap \sigma = \emptyset \mbox{ and } \tau \cup \sigma \in \Delta \}.$$
The link is itself a simplicial complex.

\item {\index{homology}\trmbitx Homology:}\quad
Associated with any simplicial complex $X$, abelian group $G$, and integer $i \ge 0$,  $H_i(X,G)$ denotes the $i$th homology group of $X$ with coefficients in $G$. If $k$ is a field, then $H_i(X,k)$ is a vector space over  $k$.

Homology is defined as ``cycles modulo boundaries''. Homology is invariant under homotopy deformations.  %The rank of $H_0(X,G)$ is equal to the number of connected components of $X$. Roughly speaking, for $i \ge 1$ $H_i(X,G)$ measures the structure of $i$-dimensional holes in $X$ as seen through the lens of the group $G$.

\item {\index{betti}\trmbitx Betti numbers:}\quad If one considers homology with coefficients in $\R$, then $H_i(X, \R)$ is a real vector space. The Betti numbers $\beta_i$ are defined by
$\beta_i = \dim \, H_i \left( X, \R \right).$ The $0$th Betti number $\beta_0$ counts the number of connected components of $X$, and in general the $i$th Betti number is said to count the number of $i$-dimensional holes in $X$.
\end{gllist}

\C{The random ${\trmitxiv 2}$-complex}
\noindent
Random hypergraphs have been well studied, but if we wish to study such objects topologically then random simplicial complexes is probably a more natural point of view.

Linial and Meshulam introduced the topological study of the random $2$-complex $Y(n,p)$ in \cite{LM06}. This model of random simplicial complex has $n$ vertices, $n \choose 2$ edges, and each of the ${n \choose 3}$ possible $2$-dimensional faces is included independently with probability $p$.

The random $2$-complex is perhaps the most natural $2$-dimensional analogue of $G(n,p)$. For example, the link of every vertex in $Y(n,p)$ has the same distribution as $G(n-1,p)$.

Several topological thresholds for $Y(n,p)$ discussed in the next section are described in Table \ref{tab:ynp}.

\vspace{-1pc}

\begin{center}
\begin{table}[htbp]
\viiipt
\baselineskip=11pt
\renewcommand{\arraystretch}{.9}
\Table{Topological thresholds for the random $2$-complex $Y=Y(n,p)$. \label{tab:ynp} c.f.\ in the TIGHT column means that the bound is best possible up to a constant factor.}
{\begin{tabular}{| p{4.5cm} | l | l | c | c |}
    \hline
\rule[-4pt]{0pt}{13pt}{\hvbviii PROPERTY}
    & {\hvbviii THRESHOLD}
    & {\hvbviii SHARP}
    & {\hvbviii TIGHT}
    & {\hvbviii SOURCE}
    \\ \hline
%\rule[0pt]{0pt}{9pt}
$Y$ is not $1$-collapsible
    & $ 2.455 / n$
    & one-sided
    & yes
    & \cite{CCFK12, ALLM13, AL13}
    \\
$H_2(Y, \R) \neq 0$
	&$ 2.753 / n$
    & one-sided
    & yes
    & \cite{Kozlov10,LP14, AL15}
    \\
$\text{cdim } \pi_1(Y) = 2$
	&$ \Theta(1 / n)$
    & ?
    & c.f.
    & \cite{CF15,Newman16}
    \\
$Y$ is not embeddable in $\R^4$
	&$ \Theta(1 / n)$
    & ?
    & c.f.
    & \cite{Wagner11}
    \\

$Y$ is pure $2$-dimensional
    & $ 2 \log n  /n$
    & sharp
    & yes
    & \cite{LM06}
    \\
$H_1(Y, \Z / \ell \Z) = 0$
    & $ 2 \log n  /n$
    & sharp
    & yes
    & \cite{LM06, MW09}
    \\
$H_1(Y, \R) = 0$
    & $ 2 \log n  /n$
    & sharp
    & yes
    & \cite{LM06,HKP12}
    \\
$\pi_1(Y)$ has property (T)
    & $ 2 \log n  /n$
    & sharp
    & yes
    & \cite{HKP12}
    \\
$H_1(Y, \Z) = 0$
    & $ O ( \log n / n )$
    & sharp
    & c.f.
    & \cite{HKP13}
    \\
$\text{cdim } \pi_1(Y)= \infty$
    & $1 / n^{3/5}$
    & coarse
    & yes
    & \cite{CF13}
    \\
$Y$ contains arbitrary subdivisions
	& $ \theta \left( 1/ \sqrt{n} \right)$
	& ?
	& c.f.
	& \cite{GW14}
	\\
$\pi_1(Y) = 0$
    & $O(1 / \sqrt{n})$
    & ?
    & no
    & \cite{BHK11,GW14, KPS16}
    \\
    \hline
\end{tabular}}
\renewcommand{\arraystretch}{1}
\xpt
\baselineskip=12pt
\end{table}
\end{center}

\vspace{-3.5pc}

\begin{table}[htbp]
\viiipt
\baselineskip=11pt
\renewcommand{\arraystretch}{.9}
\Table{Topological thresholds for $Y_d(n,p)$, $d > 2$. Definitions for the constants $ c_ d$ and $c_d^* $ are given in Section 2.3.2.
}
{\begin{tabular}{| l | l | l | c | c |}
    \hline
\rule[-4pt]{0pt}{13pt}{\hvbviii PROPERTY}
    & {\hvbviii THRESHOLD}
    & {\hvbviii SHARP}
    & {\hvbviii TIGHT}
    & {\hvbviii SOURCE}
    \\ \hline
%\rule[0pt]{0pt}{9pt}
$Y$ is not $(d-1)$-collapsible
    & $ c_d / n$
    & one-sided
    & yes
    & \cite{ALLM13,AL13}
    \\
     $H_d(Y, \R) \neq 0$
	&$ c_d^*/ n$
    & one-sided
    & yes
    & \cite{Kozlov10,LP14,AL15}
    \\
$Y$ is not embeddable in $\R^{2d}$
	&$ \Theta(1/ n)$
    & ?
    & c.f.
    & \cite{Wagner11}
    \\
$Y$ is pure $d$-dimensional
    & $ d \log n  /n$
    & sharp
    & yes
    & \cite{MW09}
    \\
$H_{d-1}(Y, \Z / \ell \Z) = 0$
    & $ d \log n  /n$
    & sharp
    & yes
    & \cite{MW09,HKP12}
    \\
$H_{d-1}(Y, \R) = 0$
    & $ d \log n  /n$
    & sharp
    & yes
    & \cite{HKP12}
    \\

$H_{d-1}(Y, \Z) = 0$
    & $ O( \log n / n)$
    & sharp
    & c.f.
    & \cite{HKP13}
    \\

$\pi_{d-1}(Y) = 0$
    & $ O( \log n / n)$
    & sharp
    & c.f.
    & \cite{HKP13}
    \\
    \hline
\end{tabular}}
\renewcommand{\arraystretch}{1}
\xpt
\baselineskip=12pt
\end{table}

\C{The random ${\trmitxiv d}$-complex}
\noindent
The natural generalization to $d$-dimensional model was introduced by Meshulam and Wallach in \cite{MW09}. For the random $d$-complex $Y_d(n,p)$,  contains the complete $(d-1)$-skeleton of a simplex on $n$ vertices, and every $d$-dimensional face appears independently with probability $p$. Some of the topological subtlety of the random $2$-dimensional model collapses in higher dimensions: for $d \ge 3$, the complexes are $d-2$-connected, and in particular simply connected. By the Hurewicz theorem, $\pi_{d-1}(Y)$ is isomorphic to $H_{d-1}(Y, \Z)$, so these groups have the same vanishing threshold.

\C{The random clique complex}
\noindent
Another analogue of $G(n,p)$ in higher dimensions was introduced  in \cite{Kahle09}. The random clique complex $X(n,p)$ is the {\it clique complex} of $G(n,p)$.
%A clique complex
It %CDT
is the maximal simplicial complex compatible with a given graph. In other words, the faces of the clique complex $X(H)$ correspond to complete subgraphs of the graph $H$.

The random clique complex asymptotically puts a measure over a wide range of topologies. Indeed, every simplicial complex is homeomorphic to a clique complex, e.g.\ by barycentric subdivision.

There are several comparisons of this model to the random $d$-complex, but some important contrasts as well. One contrast to $Y_d(n,p)$ is that for every $k \ge 1$, $X(n,p)$ has not one but two phase transitions for $k$th homology, one where homology appears and one where it vanishes. In particular, higher homology is not monotone with respect to $p$. However, there are still comparisons to $Y_d(n,p)$---the appearance of homology $H_k(X)$ is analogous to the birth of top homology $H_{d}(Y)$. Similarly, the vanishing threshold for $H_k(X(n,p))$ is analogous to vanishing of $H_{d-1}(Y)$.

\vspace{-0.5pc}

\begin{table}[htb]
\viiipt
\baselineskip=11pt
\renewcommand{\arraystretch}{.9}
\Table{Topological thresholds for $X(n,p)$}
{\begin{tabular}{| l | l | l | c | c |}
    \hline
\rule[-4pt]{0pt}{13pt}{\hvbviii PROPERTY}
    & {\hvbviii THRESHOLD}
    & {\hvbviii SHARP}
    & {\hvbviii TIGHT}
    & {\hvbviii SOURCE}
    \\  \hline 
%\rule[0pt]{0pt}{9pt}

$\pi_1(X) = 0$
    & $p = 1 / n^{1/3}$
    & ?
    & no
    & \cite{Kahle09,Babson12,CFH15}
    \\
%
%$X$ is not $k-1$-collapsible
%    & $p = c / n^{1/k}$
%    & one-sided?
%    & no
%    & \cite{Kahle09}
%    \\

$H_k(X, \R) \neq 0$
	&$ p =c / n^{1/k} $
    & one-sided
    & c.f.
    & \cite{Kahle09}
    \\
$X^k$ is pure $k$-dimensional
    & $p = \left(  \frac{(k/2 + 1 )\log n }{n} \right)^{1/(k+1)}$
    & yes
    & yes
    & \cite{Kahle09}
    \\
$H_{k}(Y, \R) = 0$
    & $p = \left(  \frac{(k/2 + 1 )\log n }{n} \right)^{1/(k+1)}$
    & yes
    & yes
    & \cite{Kahle14a}    
    \\
    [1ex]
    \hline
\end{tabular}}
\renewcommand{\arraystretch}{1}
\xpt
\baselineskip=12pt
\end{table}

\vspace{-1.5pc}

\C{The multi-parameter model}
\noindent
There is a natural multi-parameter model which generalizes all of the models discussed so far. For every every $i = 1, 2, \dots$ let $p_i: \mathbb{N} \to [0,1]$. Then define the multiparameter random complex $X(n; p_1, p_2, \dots)$ as follows. Start with $n$ vertices. Insert every edge with probability $p_1$. Conditioned on the presence of all three boundary edges, insert a $2$-face with probability $p_2$, etc.

The random $d$-complex $Y_d(n,p)$ with complete $(d-1)$-skeleton is equivalent to the case $p_1 = p_2 = \dots = p_{d-1}=1$, $p= p_d$, and $p_{d+1}=p_{d+2} = \dots = 0$. The random clique complex $X(n,p)$ corresponds to $p_2 = p_3 = \dots = 1$, and $p = p_1$. This multi-parameter model is first studied by Costa and Farber in \cite{CF14}.

\B{RANDOM GEOMETRIC MODELS}
\noindent
The \emph{random geometric graph} $\G(n,r)$ is a flexible model, defined as follows. Consider a probability distribution on $\R^d$ with a bounded, measurable, density function $f : \R^d \to \R$. Then one chooses $n$ points independently and identically distributed (i.i.d.) according to this distribution. The $n$ points are the vertices of the graph, and two vertices are adjacent if they are within distance $r$.
Usually $r = r(n) $ and $n \to \infty$. %A closely related model is to consider a Poisson point process.
The standard reference for random geometric graphs is Penrose's monograph \cite{Penrose}.

There are at least two commonly studied ways to build a simplicial complex on a geometric graph. The first is the Vietoris--Rips complex, which is the same construction as the clique complex above---one fills in all possible faces, i.e.\ the faces of the Vietoris--Rips complex correspond to the cliques of the graph. The second is the \v{C}ech complex, where one considers the higher intersections of the balls of radius $r/2$. This leads to two natural models for random geometric complexes $VR(n,r)$ and $C(n,r)$.

\vspace{-1pc}

\begin{table}[htb]
\viiipt
\baselineskip=11pt
\renewcommand{\arraystretch}{.9}
\Table{Topological thresholds for $C(n,r)$}
{\begin{tabular}{| l | l | l | c | c |}
    \hline
\rule[-4pt]{0pt}{13pt}{\hvbviii PROPERTY}
    & {\hvbviii THRESHOLD}
    & {\hvbviii SHARP}
    & {\hvbviii TIGHT}
    & {\hvbviii SOURCE}
    \\ \hline
%\rule[0pt]{0pt}{9pt}

$H_k(X) \neq 0$
	&$ nr^d =n^{-(k+2)/(k+1)} $
    & coarse
    & yes
    & \cite{Kahle11}
    \\

$H_{k}(Y) = 0$
    & $nr^d = \log n + \theta( \log \log n )$
    & sharp
    & essentially
    & \cite{BW15}
    \\
    \hline
\end{tabular}}
\renewcommand{\arraystretch}{1}
\xpt
\baselineskip=12pt
\end{table}

\vspace{-2.5pc}

\begin{table}[htb]
\viiipt
\baselineskip=11pt
\renewcommand{\arraystretch}{.9}
\Table{Topological thresholds for $VR(n,r)$}
{\begin{tabular}{| l | l | l | c | c |}
    \hline
\rule[-4pt]{0pt}{13pt}{\hvbviii PROPERTY}
    & {\hvbviii THRESHOLD}
    & {\hvbviii SHARP}
    & {\hvbviii TIGHT}
    & {\hvbviii SOURCE}
    \\ \hline
%\rule[0pt]{0pt}{9pt}

$H_k(X) \neq 0$
	&$ nr^d =n^{-(2k+2)/(2k+1)} $
    & coarse
    & yes
    & \cite{Kahle11}
    \\

$H_{k}(Y) = 0$
    & $nr^d = \theta( \log n) $
    & sharp
    & c.f.
    & \cite{Kahle11}
    \\
    \hline
\end{tabular}}
\renewcommand{\arraystretch}{1}
\xpt
\baselineskip=12pt
\end{table}

\vspace{-1.5pc}

\A{HIGH DIMENSIONAL EXPANDERS}

\vspace{-0.5pc}
\Bnn{GLOSSARY}

\begin{gllist}

\item {\index{Cheeger}\trmbitx Cheeger number:}\quad  The normalized Cheeger number of a graph $h(G)$ with vertex set $V$ is defined by
$$ h(G) = \min_{\emptyset \subsetneq A \subsetneq V} \frac{\# E( A, \bar{A}) }{\min \{ \vol  A,\vol \bar{A} \}},$$
where $$\vol A = \sum_{v \in A} \deg(v)$$
and $\bar{A}$ is the complement of $A$ in $V$.

\item {\index{Laplace}\trmbitx Laplacian:}\quad  For a connected graph $H$, the normalized graph Laplacian $L=L[H]$ is defined by
$$ L = I - D^{-1/2} A D^{-1/2}.$$
Here $A$ is the adjacency matrix, and $D$ is the diagonal matrix with vertex degrees along the diagonal.

\item {\index{Gap}\trmbitx Spectral gap:} The eigenvalues of the normalized graph Laplacian of a connected graph satisfy
$$ 0 = \lambda_1 < \lambda_2 \le \dots \le \lambda_n \le 2.$$
The smallest positive eigenvalue $\lambda_2[H]$ is of particular importance, and is sometimes called the spectral gap of $H$.

\item {\index{Expander}\trmbitx Expander family:} Let $\{ G_i \}$ be an infinite sequence of graphs where the number of vertices tends to infinity. We say that $\{ G_i \}$ is an expander family if
$$  \liminf \lambda_2 [G_i] > 0.$$

\end{gllist}

Expander graphs are of fundamental importance for their applications in computer science and mathematics \cite{HLW06}. It is natural to seek their various higher-dimensional generalizations. See \cite{Lubotzky14} for a survey of recent progress on higher-dimensional expanders, particularly Ramanujan complexes which generalize Ramanujan graphs.

Gromov suggested that one property that higher-dimensional expanders should have is geometric or topological overlap. A sequence of $d$-dimensional simplicial complexes $\Delta_1$, $\Delta_2$, \dots, is said to have the \emph{geometric overlap property} if for every geometric map (affine-linear on each face) $f: \Delta_i \to \R^d$, there exists a point $p \in \R^d$ such that $f^{-1}(p)$ intersects the interior of a constant fraction of the $d$-dimensional faces.  The sequence is said to have the stronger \emph{topological overlap property} if this holds even for continuous maps.

One way to define higher-dimensional expander is via \emph{coboundary expansion}, which generalizes the Cheeger number of a graph. Following Linial and Meshulam's coisoperimetric ideas, Dotterrer and Kahle \cite{DK12} pointed out that $d$-dimensional random simplicial complexes are coboundary expanders. By a theorem of Gromov \cite{Gromov09,Gromov10} (see the note \cite{DKW15} for a self contained proof of Gromov's theorem), this implies that random complexes have the topological overlap property.

Lubotzky and Meshulam introduced a new model of random $2$-complex, based on random Latin squares, in \cite{LM15}. The main result is the existence of coboundary expanders with bounded edge-degree, answering a question asked implicitly in \cite{Gromov10} and explicitly in \cite{DK12}.

Another way to define higher-dimensional expanders is via \emph{spectral gap} of various Laplacian operators. Hoffman, Kahle, and Paquette studied spectral gap of random graphs \cite{HKP12}, and applied Garland's method to prove homology-vanishing theorems. Gundert and Wagner extened this to study higher-order spectral gaps of these complexes \cite{GW15}. Parzanchevski, Rosenthal, and Tessler showed that this implies the geometric overlap property \cite{PRT15}.

%Random simplicial complexes provide higher-dimensional generalizations of random graphs, and so we see many generalizations of random graph theory to higher dimensions. Random complexes have interesting expander-like properties (Cheeger constant, spectral, geometric overlap) . \cite{Lubotzky14,DK12,Gromov09,Gromov10,GW12}
%

\clearpage

\A{PHASE TRANSITIONS}
\noindent
There has been a lot of interest in identifying thresholds for various topological properties, such as vanishing of homology. As some parameter varies,  the topology passes a \emph{phase transition} where some property suddenly emerges. In this section we review a few of the most well-studied topological phase transitions.

\B{HOMOLOGY-VANISHING THEOREMS}

%\C{The Erd\H{o}s--R\'enyi theorem}
\noindent
The following theorem describing the connectivity threshold for the random graph $G(n,p)$ is the archetypal homology-vanishing theorem.

\theorem{Erd\H{o}s--R\'enyi theorem \, {\cmrx\cite{ER59}}}
{If $$ p \ge \frac{ \log n + \omega(1) }{n} $$ then w.h.p.\ $G(n,p)$ is connected, and
if $$ p \le \frac{ \log n - \omega(1) }{n} $$ then w.h.p.\ $G(n,p)$ is disconnected.
Here $\omega(1)$ is any function so that $\omega(1) \to \infty$ as $n \to \infty$.}

We say that it is a homology-vanishing theorem because path connectivity of a topological space $X$ is equivalent to $\widetilde{H}_0(X, G) =0$ with any coefficient group $G$.

It is also a cohomology-vanishing theorem, since $\widetilde{H}^0(X, G) =0$ is also equivalent to path-connectivity. In many ways it is better to think of it as a cohomology theorem, since the standard proof, for example in Chapter 10 of \cite{Bollobas}, is really a cohomological one.
This perspective helps when understanding the proof of the Linial--Meshulam theorem.

%\C{The Linial--Meshulam theorem}
%\noindent
The following cohomological analogue of the Erd\H{o}s--R\'enyi theorem was the first nontrivial result for the topology of random simplicial complexes.

\theorem{Linial--Meshulam theorem \, {\cmrx\cite{LM06}\label{thm:LM}}}
{Let $Y = Y(n,p)$. If
$$ p \ge \frac{2 \log n + \omega(1) }{n} $$
then w.h.p.\ $H_1(Y, \Z / 2\Z) =0$, and if
$$ p \le \frac{ 2 \log n - \omega(1) }{n} $$
then w.h.p.\ $H_1(Y, \Z / 2\Z)  \neq 0$.
}

\smallskip
One of the main tools introduced in \cite{LM06} is a new co-isoperimetric inequality for the simplex, which was discovered independently by Gromov. These co-isoperimetric inequalities were combined by Linial and Meshulam with intricate cocycle-counting combinatorics to get a sharp threshold.

See \cite{DKW15} for a comparison of various definitions co-isoperimetry, and a clean statement and self-contained proof of Gromov's theorem.

Theorem \ref{thm:LM} was generalized further by Meshulam and Wallach.

\theorem{\cmrx\cite{MW09} \label{thm:MW}}
{Fix $d \ge 1$, and let $Y = Y_d(n,p)$. Let $G$ be any finite abelian group.
If $$ p \ge \frac{d \log{n} + \omega(1) }{n} $$ then w.h.p.\ $H_{d-1}(Y, G) =0$, and
if $$ p \le \frac{ d \log{n} - \omega(1) }{n} $$ then w.h.p.\ $H_{d-1}(Y, G)  \neq 0$.}

Theorem \ref{thm:MW} generalizes Theorem \ref{thm:LM} in two ways: by letting the dimension $d \ge 2$ be arbitrary, and also by letting coefficients be in an arbitrary finite abelian group $G$.

\C{Spectral gaps and Garland's method}
\noindent
There is another approach to homology-vanishing theorems for simplicial complexes, via Garland's method \cite{Garland73}.
The following refinement of Garland's theorem is due to Ballman and \'Swiatkowski.

\theorem{\cmrx\cite{Garland73,BS97}}
{If $\Delta$ is a finite, pure $d$-dimensional, simplicial complex, such that
$$\lambda_2 [\link_{\Delta} (\sigma)] > 1-  \frac{1}{d}$$ for every $(d-2)$-dimensional face $\sigma \in \Delta$,
then $H_{d-1}(\Delta, \R) =0$.}
This leads to a new proof of Theorem \ref{thm:MW}, at least over a field of characteristic zero.

\theorem{\cmrx\cite{HKP12}}
{Fix $d \ge 1$, and let $Y = Y_d(n,p)$.
If $$ p \ge \frac{d \log{n} + \omega(1) }{n} $$ then w.h.p.\ $H_{d-1}(Y, \R) =0$, and
if $$ p \le \frac{ d \log{n} - \omega(1) }{n} $$ then w.h.p.\ $H_{d-1}(Y, \R)  \neq 0$.
Here $\omega(1)$ is any function that tends to infinity as $n \to \infty$.}

This is slightly weaker than the Meshulam--Wallach theorem topologically speaking, since $H_i(Y, G) = 0$ for any  finite group $G$ implies that $H_i(Y, \R) = 0$ by the universal coefficient theorem, but generally the converse is false. However, the proof via Garland's method avoids some of the combinatorial complications of cocycle counting. %The method also provides sharp probabilistic information such as Poisson distribution of Betti numbers in the critical window, hitting-time results, etc.

Garland's method also provides proofs of theorems which have so far eluded other methods. For example, we have the following homology-vanishing threshold in the random clique complex model. Note that $k=0$ again corresponds to the Erd\H{o}s--R\'enyi theorem.

\theorem{\cmrx\cite{Kahle14a} \label{thm:sharp}}
{Fix $k \ge 1$ and let $X = X(n,p)$. Let $\omega(1)$ denote a function that tends to $\infty$ arbitrarily slowly.
If $$ p \ge \left( \frac{ \left( \frac{k}{2}+1 \right) \log{n} + \left( \frac{k}{2}  \right) \log\log{n} + \omega(1)}{n} \right)^{1 / (k+1)}$$
then w.h.p.\ $H_k(X,\R) = 0$, and if
$$\frac{1}{n^k} \le  p \le \left( \frac{ \left( \frac{k}{2}+1 \right) \log{n} + \left( \frac{k}{2}  \right) \log\log{n} - \omega(1)}{n} \right)^{1 / (k+1)}$$
then w.h.p.\ $H_k(X,\R) \neq 0$.}

\smallskip
It may be that this theorem holds with $\R$ coefficients replaced by a finite group $G$ or even with $\Z$, but for the most part this remains an open problem. The only other case that seems to be known is the case $k=1$ and $G  = \Z / 2$ by DeMarco, Hamm, and Kahn \cite{DHK13}, where a similarly sharp threshold is obtained.

\smallskip

The applications of Garland's method depends on new results on the spectral gap of random graphs.

\theorem{\cmrx\cite{HKP12} \label{thm:gap}}
{Fix $k \ge 0$. Let $\lambda_1, \lambda_2, \dots$ denote the eigenvalues of the normalized graph Laplacian of the random graph $G(n,p)$.
If
$$p \ge  \frac{ (k+1) \log n + \omega(1)}{n}$$
then
$$1- \sqrt{\frac{C}{np}} \le \lambda_2 \le \dots \le \lambda_n \le 1+ \sqrt{\frac{C}{np}}$$ with probability at least $1 - o \left(n^{-k} \right)$. Here $C >0$ is a universal constant.}

Theorem \ref{thm:sharp}, combined with some earlier results \cite{Kahle09}, has the following corollary.
\theorem{\cmrx\cite{Kahle14a}}
{Let $k \ge 3$ and $\epsilon > 0$ be fixed.
If $$ \left( \frac{ \left(C_k + \epsilon \right) \log n }{n} \right)^{1/k}  \le p \le \frac{1}{n^{1/(k+1)+ \epsilon}},$$
where $C_3 = 3$ and $C_k =   k/2+1$ for $k > 3$, then w.h.p.\ $X$ is rationally homotopy equivalent to a bouquet of $k$-dimensional spheres.\\}

The main remaining conjecture for the topology of random clique complexes is that these rational homotopy equivalences are actually homotopy equivalences.
\conj{The bouquet-of-spheres conjecture.}
{Let $k \ge 3$ and $\epsilon > 0$ be fixed.  If
$$ \frac{n^{\epsilon} }{n^{1/k}} \le p \le \frac{n^{-\epsilon}}{n^{1/(k+1)}}$$
then w.h.p. $X$ is homotopy equivalent to a bouquet of $k$-spheres.}

Given earlier results, this is equivalent to showing that $H_k(X, \Z)$ is torsion free. So far, integer homology for $X(n,p)$ is not very well understood. Some progress has been made for $Y_d(n,p)$, described in the following.

\C{Integer homology}
\noindent
Unfortunately, neither method discussed above (the cocycle-counting methods pioneered by Linial and Meshulam or the spectral methods of Garland), seems to handle integral homology. There is a slight subtlety here---if one knows for some simplicial complex $\Sigma$ that $H_i(\Sigma, G) = 0$ for every finite abelian group $G$ then $H_i(\Sigma, \Z) = 0$ by the universal coefficient theorem. See, for example, Chapter 2 of Hatcher \cite{Hatcher}.

So it might seem that the Theorem \ref{thm:MW} will also handle $\Z$ coefficients, but the proof uses cocycle counting methods which require $G$ to be fixed, or at least for the order of the coefficient group $|G|$ to be growing sufficiently slowly. Cocycle counting does not seem to work, for example, when $|G|$ is growing exponentially fast. The following gives an upper bound on the vanishing threshold for integer homology.

\theorem{\cmrx\cite{HKP12}}
{Fix $d \ge 2$, and let $Y = Y_d(n,p)$.
If $$ p \ge \frac{80d \log{n} }{n},$$ then w.h.p.\ $H_{d-1}(Y, \Z) =0$.}

The author suspects that the true threshold for homology with $\Z$ coefficients is the same as for field coefficients: $d \log n / n$.
%The technique in \cite{HKP12} is elementary but seems to be a new approach to homology-vanishing theorems, circumventing some of the technical obstacles inherent in cocycle counting and spectral gap methods.

\conj{A sharp threshold for $\Z$ homology.}
{If $$ p \ge \frac{d \log{n} + \omega(1) }{n} $$ then w.h.p.\ $H_{d-1}(Y, \Z) =0$, and
if $$ p \le \frac{ d \log{n} - \omega(1) }{n} $$ then w.h.p.\ $H_{d-1}(Y, \Z)  \neq 0$.}

\B{THE BIRTH OF CYCLES AND COLLAPSIBILITY}

\vspace{-1pc}

\C{${\trmitxiv G(n,p)}$ in the ${\trmitxiv p=1/n}$ regime}
\noindent
There is a remarkable phase transition in structure of the random graph $G(n,p)$ at the threshold $p=1/n$. A ``giant'' component, on a constant fraction of the vertices, suddenly emerges. This is considered an analogue of percolation on an infinite lattice, where an infinite component appears with probability $1$.

\theorem{\cmrx\cite{ER59} \label{thm:giant}}
{Let $p = c / n$ for some $c > 0$ fixed, and $G = G(n,p)$.
\begin{itemize}
\item If $c < 1$ then w.h.p.\ all components are of order $O( \log n)$.
\item If $c > 1$ then w.h.p.\ there is a unique giant component, of order $\Omega(n)$.
\end{itemize}
}

An overview of this remarkable phase transition can be found in Chapter 11 of Alon and Spencer \cite{AlonSpencer}.

In random graphs, the appearance of cycles with high probability has the same threshold $1/n$.

\theorem{\cmrx\cite{Pittel88} \label{thm:cycles}}
{Suppose $p = c / n$ where $c > 0 $ is constant.
\begin{itemize}
\item If $c \ge 1$ then w.h.p. $G$ contains at least one cycle, i.e.
 $$\prob \left[ H_1(G) \neq 0 \right] \to 1.$$
\item If $c < 1$ then
$$\prob \left[ H_1(G) = 0 \right]  \to \sqrt{1-c} \, \exp(c/2 + c^2 / 4).$$
\end{itemize}
}

The analogy in higher dimensions is only just beginning to be understood.

\C{The birth of cycles}
\noindent
Kozlov first studied the vanishing threshold for top homology in  \cite{Kozlov10}.

\theorem{\cmrx\cite{Kozlov10}}
{Let $Y = Y_d(n,p)$, and $G$ be any abelian group.
\begin{itemize}
\item[{\rm (1)}] If $p= o(1/n)$ then w.h.p.\ $H_d(Y, G) = 0$.
\item[{\rm (2)}] If $p = \omega(1 / n)$ then w.h.p.\ $H_d(Y, G) \neq 0$.
\end{itemize}
}

Part (1) of this theorem cannot be improved. Indeed, let $S$ be the number of subcomplexes isomorphic to the boundary of a $(d+1)$-dimensional simplex. If $p = c / n$ for some constant $c > 0$, then
$$\expect[S] \to c^{d+2} / (d+2)!,$$
as $n \to \infty$. Moreover, $S$ converges in law to a Poisson distribution with this mean in the limit, so
$$ \prob[ H_d(Y, G)  \neq 0] \ge \prob[S \neq 0] \to 1 - \exp (-c^{d+2} / (d+2)! ).$$
In particular, for $p = c/n$ and $c >0$, $\prob[ H_d(Y, G) \neq 0]$ is bounded away from zero.

On the other hand, part (2) can be improved. Indeed straightforward computation shows that if $p \ge c / n$ and $c > d+1$ then w.h.p.\ the number of $d$-dimensional faces is greater than the number of $(d-1)$-dimensional faces. Simply by dimensional considerations, we conclude that $H_d(Y, G) \neq 0$.

This can improved more though. Aronshtam and Linial found the best possible constant factor $c_d^*$, defined for $d \ge 2$ as follows.

Let $x \in (0,1)$ be the unique root to the equation
$$(d+1)(1-x) + (1 + dx) \log x = 0,$$
and then set
$$c_d^* = \frac{ -\log x}{ (1 - x)^d}.$$
\theorem{\cmrx\cite{AL15}}
{Let $Y =Y_d(n,p)$. If $p \ge c / n$ where $c > c_d^*$, then w.h.p.\ $H_d(Y, G) \neq 0$.}

\smallskip
In the other direction, Linial and Peled showed that this result is tight, at least in the case of $\R$ coefficients.

\theorem{\cmrx\cite{LP14}}
{If $p \le c / n$ where $c < c_d^*$ then w.h.p.\ $H_d(Y, G)$ is generated by simplex boundaries. So
$$\prob \left[ H_d(Y, \R) = 0 \right] \to \exp (-c^{d+2} / (d+2)! ).$$
}

Linial and Peled also showed the birth of a \emph{giant (homological) shadow} at the same point. This is introduced and defined in \cite{LP14}, and it is discussed there as a higher-dimensional analogue of the birth of the giant component in $G(n,p)$.

\C{The threshold for $d$-collapsibility}
\noindent
In a $d$-dimensional simplicial complex, an \emph{elementary collapse} is an operation that deletes a pair of faces $(\sigma, \tau)$ such  that
\begin{itemize}
\item $\tau$ is a $d$-dimensional face,
\item $\sigma$ is a $d-1$-dimensional face contained in $\tau$, and
\item $\sigma$ is not contained in any other $d$-dimensional faces.
\end{itemize}

An elementary collapse results in a homotopy equivalent simplicial complex.

\medskip
If a simplicial complex can be reduced to a $d-1$-dimensional complex by a series of elementary collapses, we say that it is \emph{d-collapsible}.

For a graph, $1$-collapsible is equivalent to being a forest. In other words, a graph $G$ is $1$-collapsible if and only if $H_1(G)=0$.
This homological criterion does not hold in higher dimensions. In fact, somewhat surprisingly, $d$-collapsibility and $H_d \neq 0$ have distinct thresholds for random complexes.

Let $d \ge 2$ and set
$$g_d(x) = (d+1) (x+1)e^{-x} + x (1-e^{-x})^{d+1}.$$
Define $c_d$ to be the unique solution $x > 0$ of $g_d(x) = d+1$.

\theorem{\cmrx\cite{ALLM13, AL13} \label{thm:collapse}}
{Let $Y = Y_d(n,p)$.
\begin{itemize}
\item If $p \ge c / n$ where $c > c_d$ then w.h.p.\ $Y$ is not $d$-collapsible, and
\item if $p \le c / n$ where $ c < c_d$ then $Y$ is $d$-collapsible with probability bounded away from zero.
\end{itemize}
}

So again, this is a one-sided sharp threshold.
Regarding collapsibility in the random clique complex model, Malen showed in his PhD thesis that if $p \ll n^{-1 / (k+1)}$, then w.h.p.\ $X(n,p)$ is $k$-collapsible.

\C{Embeddability}
\noindent
Every $d$-dimensional simplicial complex is embeddable in $\R^{2d+1}$, but not necessarily in $\R^{2d}$. Wagner studied the threshold for non-embeddability of random $d$-complexes in $\R^{2d}$, and showed the following for $Y = Y_d(n,p)$.

\theorem{\cmrx\cite{Wagner11} \label{thm:embed}}
{There exists constants $c_1, c_2 > 0$ depending only on the dimension $d$ such that:
\begin{itemize}
\item if $p < c_1 / n$ then w.h.p.\ $Y$ is embeddable in $\R^{2d}$, and
\item if $p > c_2 / n$ then w.h.p.\ $Y$ is not embeddable in $\R^{2d}$.
\end{itemize}
}

There is a folklore conjecture that a $d$-dimensional simplicial complex on $n$ vertices embeddable in $\R^{2d}$ can have at most $O(n^d)$ faces \cite{Dey93,Kalai91}. See, for example, the discussion in the expository book chapter \cite{Wagner13} or Chapter 24 of this handbook. The $d=1$ case is equivalent to showing that a planar graph may only have linearly many edges, which follows immediately from the Euler formula, but the conjecture is open for every $d \ge 2$. Theorem \ref{thm:embed} shows that it holds generically.

\B{PHASE TRANSITIONS FOR HOMOLOGY IN RANDOM\\ GEOMETRIC COMPLEXES}

\smallskip\noindent
Penrose described sharp thresholds for connectivity of random geometric graphs, analogous to the Erd\H{o}s--R\'enyi theorem. In the case of a uniform distribution on the unit cube $[0,1]^d$ or a standard multivariate distribution, these results are tight \cite{Penrose}.

Thresholds for homology in random geometric complexes was first studied in \cite{Kahle11}. A homology vanishing threshold for random geometric complexes is obtained in \cite{Kahle11}, which is tight up to a constant factor, but recently a much sharper result was obtained by  Bobrowski and Weinberger.
\theorem{\cmrx\cite{BW15} }
{Fix $1 \ge k \ge d-1$. If $$nr^d \ge \log n + k \log \log n + \omega(1),$$ then w.h.p. $\beta_k = 0$, and if
$$nr^d \le \log n + (k-2) \log \log n - \omega(\log \log \log n),$$
then w.h.p. $\beta_k \to \infty$.}

\B{RANDOM FUNDAMENTAL GROUPS}

\vspace{-0.5pc}
\Bnn{GLOSSARY}

\begin{gllist}

\item {\index{fundamental!group}\trmbitx Fundamental group:}\quad
In a path-connected topological space $X$, choose an arbitrary base point $p$. Then the homotopy classes of loops in $X$ based at $p$, i.e.\ continuous functions $f: [0,1] \to X$ with $f(0)=f(1)=p$ may be endowed with the structure of a group, where the group operation is concatenation of two loops at double speed. This is called the fundamental group $\pi_1(X)$, and up to isomorphism it does not depend on the choice of base point $p$. If $\pi_1(X) = 0$ then $X$ is said to be simply connected. The first homology group $H_1(X, \Z)$ is isomorphic to the abelianization of $\pi_1(X)$.

\item {\index{chain}\trmbitx A chain of implications:}\quad  The following implications hold for an arbitrary simplicial complex $X$.
$$\pi_1(X) = 0 \implies H_1(X, \Z) = 0 \implies H_1(X, \Z / q \Z) = 0  \implies H_1(X, \R) = 0.$$
Here $q$ is any prime. This is a standard application of the universal coefficient theorem for homology \cite{Hatcher}.

A partial converse to one of the implications is the following. If $H_1(X, \Z / q \Z) = 0$ for every prime $q$, then $H_1(X, \Z) = 0$.

\item {\index{hyper}\trmbitx Hyperbolic group:}\quad A finitely presented group is said to be word hyperbolic if it can be equipped with a word metric satisfying certain characteristics of hyperbolic geometry \cite{Gromov87}.

\item  {\index{T}\trmbitx Kazhdan's property (T):}\quad A group $G$ is said to have property (T) if the trivial representation is an isolated point in the unitary dual equipped with the Fell topology. Equivalently, if a representation has almost invariant vectors then it has invariant vectors.

\item {\index{Gcohom}\trmbitx Group cohomology:}\quad Associated with a finitely-presented group $G$ is a contractible CW complex $EG$ on which $G$ acts freely. The quotient $BG$ is the classifying space for principle $G$ bundles. The group cohomology of $G$ is equivalent to the cohomology of $BG$.

\item {\index{cdim}\trmbitx Cohomological dimension:}\quad  The cohomological dimension of a group $G$, denoted $\mbox{cdim }G$, is the largest dimension $k$ such that $H^k ( G, R) \neq 0$ for some coefficient ring $R$.

\end{gllist}

The random fundamental group $\pi_1(Y(n,p))$ may fruitfully be compared to other models of random group studied earlier, such as Gromov's density model \cite{Ollivier05}. The techniques and flavor of the subject owes as much to geometric group theory as to combinatorics.

\C{The vanishing threshold and hyperbolicity}
\noindent
Babson, Hoffman, and Kahle showed that the vanishing threshold for simple connectivity is much larger than the homology-vanishing threshold.

\theorem{\cmrx\cite{BHK11} \label{thm:fun}}
{Let $\epsilon >  0$ be fixed and $Y = Y(n,p)$.
If $$p \ge \frac{n^{\epsilon}}{\sqrt{n}}$$ then w.h.p.\ $\pi_1(Y) = 0$,
and if $$p \le \frac{n^{ - \epsilon}}{\sqrt{n}}$$ then w.h.p.\ $\pi_1(Y)$ is a nontrivial hyperbolic group.}

Most of the work in proving Theorem \ref{thm:fun} is showing that, on the sparse side of the threshold, $\pi_1$ is hyperbolic. This in turn depends on a local-to-global principle for hyperbolicity due to Gromov \cite{Gromov87}.

Gundert and Wagner showed that it suffices to assume that $$p \ge \frac{C}{\sqrt{n}}$$ for some constant $C > 0$ to show that w.h.p.\ $\pi_1(Y) = 0$ \cite{GW14}.  Kor\'andi, Peled, and Sudakov showed that it suffices to take $C=1/2$ \cite{KPS16}.

The author suspects that there is a sharp threshold for simple connectivity at $C / \sqrt{n}$ for some $C>0$.

\conj{A sharp vanishing threshold for $\pi_1(Y)$.}
{There exists some constant $C > 0$ such that if
$$p \ge \frac{C+\epsilon}{\sqrt{n}},$$
with high probability, $\pi_1(Y) = 0$; and if
$$p \le \frac{C-\epsilon}{\sqrt{n}},$$
with high probability, $\pi_1(Y) \neq 0$.
}

\C{Kazhdan's property (T)}
\noindent
One of the most important properties studied in geometric group theory is property (T).  Loosely speaking, a group is (T) if it does not have many unitary representations. Property (T) is also closely related to the study of expander graphs. For a comprehensive overview of the subject, see the monograph \cite{BdlHV}.

Inspired by Garland's method, \.Zuk gave a spectral condition sufficient to imply (T). Hoffman, Kahle, and Paquette applied \.Zuk's condition, together with Theorem \ref{thm:gap} to show that the threshold for $\pi_1(Y)$ to be (T) coincides with the Linial--Meshulam homology-vanishing threshold.\\

\theorem{\cmrx\cite{HKP12} \label{thm:T}}
{Let $Y = Y(n,p)$.
\begin{itemize}
\item If $$ p \ge \frac{2\log n + \omega(1) }{n} $$ then w.h.p.\ $\pi_1(Y)$ is (T), and
\item if $$ p \le \frac{ 2\log n - \omega(1) }{n} $$ then w.h.p.\ $\pi_1(Y)$ is not (T).
\end{itemize}
}

\C{Cohomological dimension}
\noindent
Costa and Farber~\cite{CCFK12} studied the cohomological dimension of the random fundamental group in \cite{CF13}. Their main findings are that there are regimes when the cohomological dimension is $1$, $2$, and $\infty$, before the collapse of the group at $p = 1 / \sqrt{n}$.

\theorem{\cmrx\cite{HKP12} \label{thm:cdim}}
{Let $Y = Y(n,p)$.
\begin{itemize}
\item If $$p \ll \frac{1}{n}$$ then w.h.p.\ $\mbox{cdim } \pi_1(Y)=1$ \cite{CCFK12}.
\item If $$ \frac{3 }{n} \le p \ll n^{-3/5}$$ then w.h.p.\ $\mbox{cdim } \pi_1(Y)=2$ \cite{CF15} .
\item if $$n^{-3/5} \ll p \le n^{-1/2 - \epsilon}$$ then w.h.p.\ $\mbox{cdim } \pi_1(Y) = \infty$ \cite{CF13}.
\end{itemize}
}
Here we use $f \ll g $ to mean $f= o(g)$, i.e.
$$\lim_{n \to \infty} f/g = 0.$$

Newman recently refined part of this picture \cite{Newman16},  showing that if $p < 2.455 / n$ w.h.p.\ $\mbox{cdim } \pi_1(Y)=1$, and if $p > 2.754 / n$ then w.h.p.\ $\mbox{cdim } \pi_1(Y)=2$. The precise constants are $c_2$ and $c_2^*$, defined in Section 23.3.2.

\C{The fundamental group of the clique complex}
\noindent
Babson showed that $p=n^{-1/3}$ is the vanishing threshold for $\pi_1( X(n,p) )$ in  \cite{Babson12}. An independent and self-contained proof, including more refined results regarding torsion and cohomological dimension, was given by Costa, Farber, and Horak in~\cite{CFH15}.

\C{Finite quotients}
\noindent
Meshulam studied finite quotients of the random fundamental group, and showed that if they exist then the index must be large---the index must tend to infinity with $n$. His technique is a version of the cocycle-counting arguments in \cite{LM06} and \cite{MW09}, for non-abelian cohomology.

\theorem{\cmrx Meshulam, \cite{Meshulam13}}
{Let $ c > 0 $ be fixed.  If $ p \ge \frac{(6 + 7c) \log n}{n} $
then w.h.p. $\pi_1(Y)$ has no finite quotients with index less than $n^c$. Moreover, if $H$ is any \emph{fixed} finite group and $ p \ge \frac{(2 + c) \log n}{n} $ then w.h.p.\ there are no nontrivial maps to $H$.}

\B{PHASE TRANSITIONS IN THE MULTI-PARAMETER MODEL}
\noindent
Applying Garland's method, Fowler described the homology-vanishing phase transition in the multi-parameter model in \cite{Fowler15}.

\theorem{\cmrx Fowler \cite{Fowler15}}
{Let $X = X(n, p_1, p_2, \dots)$ with $p_i = n^{- \alpha_i}$ and $\alpha_i \ge 0$ for all $i$. If
$$\sum_{i=1}^k \alpha_i {k \choose i} < 1,$$
then w.h.p.\ $H^{k-1}(X, \Q) = 0$.
\medskip
If $$\sum_{i=1}^k \alpha_i {k \choose i} \ge 1$$
and $$\sum_{i=1}^{k-1} \alpha_i {k-1 \choose i} < 1$$
then w.h.p. $H^{k-1}(X, \Q) \neq 0$.}

\A{BETTI NUMBERS AND PERSISTENT HOMOLOGY}

%%Rather than ask whether homology is vanishing or non-vanishing, a more refined question is to ask how large it is in the non-vanishing regime. Betti numbers is one measure of this.
%
%The $i$th Betti number of a topological space $X$ is defined by $$\beta_i (X) = \dim \, H_i ( X, \R).$$
%
%It is well known that random graphs have strong expander-like properties.
%
%For example, the normalized Cheeger number of $G(n,p)$ is bounded away from zero, once $p \gg \log n / n$.
%
%Similarly, the spectral gap of the normalized graph Laplacian of $G(n,p)$ is close to $1$ in the same regime, by the $k=0$ case of Theorem \ref{thm:gap}.
%

%\cite{GW15}

\vspace{-0.5pc}
\B{BETTI NUMBERS}

\vspace{-0.5pc}
\C{The random clique complex}
\noindent
In the random clique complex $X(n,p)$, it was noted in \cite{Kahle09} that if $$1/ n^{1/k} \ll p \ll 1 / n^{1/(k+1)},$$
 then
$$\expect \left[ \beta_k \right] = \left(1 - o(1) \right) {n \choose k+1} p^{k+1 \choose 2}.$$
A more refined estimate may be obtained along the following lines.

We have the Euler relation
$$ \chi = f_0 - f_1 + f_2 - \dots = \beta_0 - \beta_1 + \beta_2 - \dots,$$
where $f_i$ denotes the number of $i$-dimensional faces.
The expected number of $i$-dimensional faces is easy to compute---by linearity of expectation we have
$$\expect[ f_i ] = {n \choose i+1 } p ^{i+1 \choose 2}.$$
If we make the simplifying assumption that only one Betti number $\beta_i$ is nonzero, then we have
$$ | \chi | = | f_0 - f_1 + f_2 - \dots | = \beta_i.$$

\begin{figure}
\centering
\includegraphics[width=3in,resolution=600]{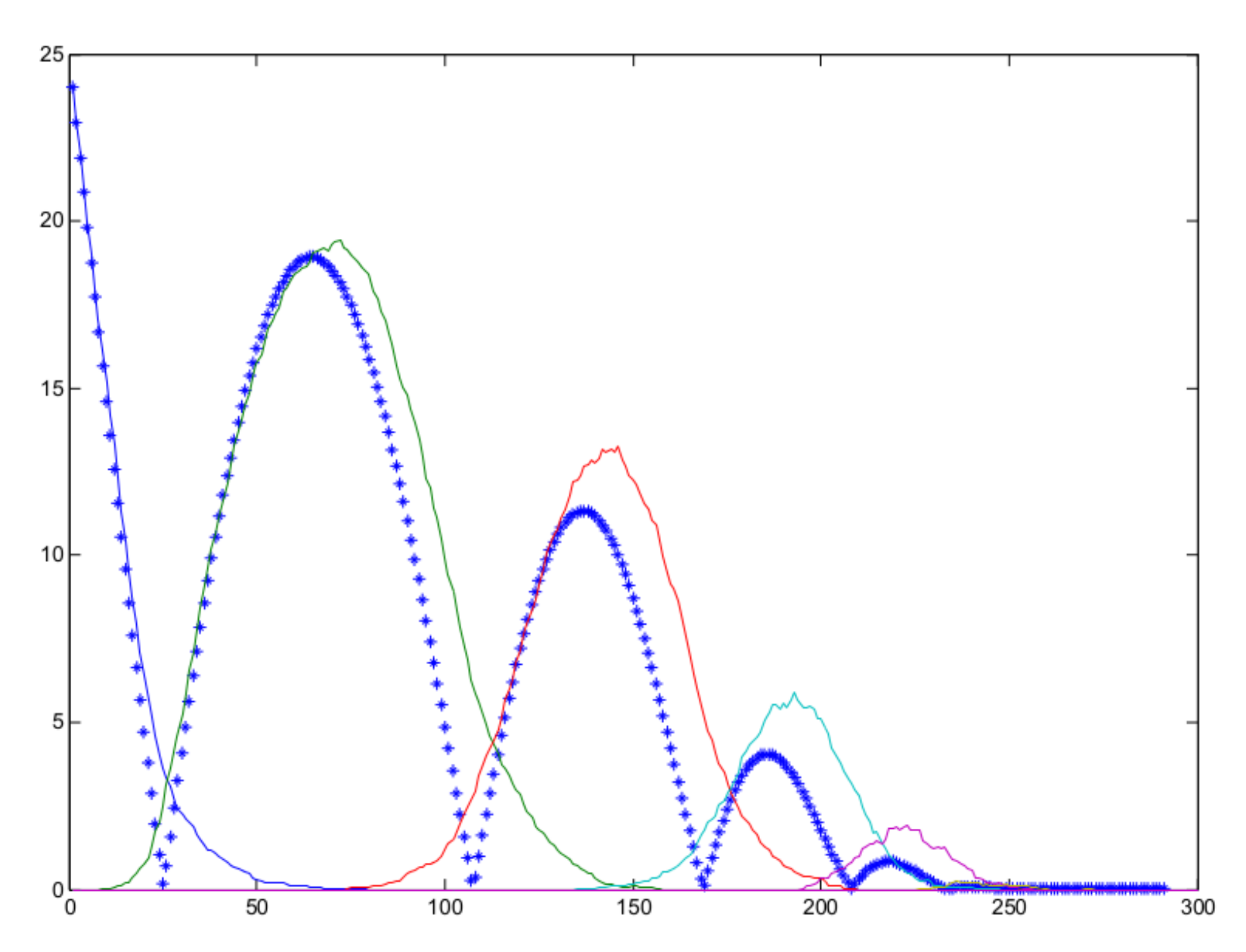}
\caption{$|\expect[ \chi ]|$ plotted in blue, against the Betti numbers of the random flag complex $X(n,p)$. Here $n=25$ and $p$ varies from $0$ to $1$. The horizontal axis is the number of edges. Computation and image courtesy of Vidit Nanda.}
\label{fig:Betti}
\end{figure}

So we obtain a plot of all of the Betti numbers by plotting the single function
\begin{align*}
 | \chi | &= | f_0 - f_1 + f_2 - \dots | \\
 & = \left| {n \choose 1} - {n \choose 2} p^1 + {n \choose 3} p^2 - {n \choose 4} p^6 + \dots \right|
\end{align*}
This seems to work well in practice. See for example Figure \ref{fig:Betti}. It is interesting that even through all of the theorems we have discussed are asymptotic as $n \to \infty$, the above heuristic gives a reasonable prediction of the shape of the Betti number curves, even for $n=25$ and $p \le 0.6$.

\C{Homological domination in the multi-parameter model}
\noindent
In \cite{CF15dom, CF15crit}, Costa and Farber show that for many choices of parameter (an open, dense subset of the set of allowable vectors of exponents) in the multi-parameter model, the homology is dominated in one degree.

\C{Random geometric complexes}
\noindent
Betti numbers of random geometric complexes were first studied by Robins in \cite{Robins06}.

Betti numbers of random geometric complexes are also studied in \cite{Kahle11}. Estimates are obtained for the Betti numbers in the subcritical regime $nr^d \to 0$. In this regime the Vietoris--Rips complex and \v{C}ech complex have small connected components (bounded in size), so all the topology is local.

For the following theorems, we assume that $n$ points are chosen i.i.d.\ according to a probability measure on $\R^d$ with a bounded measurable density function $f$. So the assumptions on the underlying probability distribution are farily mild.

The following describes the expectation of the Betti numbers of the Vietoris--Rips complex in the subcritical regime. In this regime the homology of $VR(n,r)$ is dominated by subcomplexes combinatorially isomorphic to the boundary of the cross-polytope.

\theorem{\cmrx\cite{Kahle11} \label{thm:BettiRips}}
{Fix $d \ge 2$ and $k \ge 1$. If $nr^d \to 0$ then
$$\expect \left[ \beta_k [VR(n,r)] \right] / \left( n^{2k+2}r^{d(2k+1)} \right) \to C_k,$$
as $n \to \infty$ where $D_k$ is a constant which depends on $k, d$, and the function $f$.}

The analogous story for the \v{C}ech complex is the following. Here the homology is dominated by simplex boundaries.

\theorem{\cmrx\cite{Kahle11} \label{thm:BettiCech}}
{Fix $d \ge 2$ and $1 \le k \le d-1$. If $nr^d \to 0$ then
$$\expect \left[ \beta_k [ C(n,r) ] \right] / \left( n^{k+2}r^{d(k+1)} \right) \to C_k,$$
as $n \to \infty$ where $D_k$ is a constant which depends on $k, d$, and the function $f$.}

\C{In the thermodynamic limit}
\noindent
The thermodynamic limit, or critical regime, is when $nr^d \to C$ for some constant $C > 0$.
In \cite{Kahle11}, it is shown that for every $1 \le k \le d-1$, we have $\beta_k = \Theta(n)$.
Yogeshwaran, Subag, and Adler obtained the strongest results so far for Betti numbers in the thermodynamic limit, including strong laws of large numbers \cite{YSA15}, in particular that $\beta_k / n \to C_k$.

\C{Limit theorems}
\noindent
Kahle and Meckes computed variance of the Betti numbers, and proved Poisson and normal limiting distributions for Betti numbers in the subcritical regime $r=o \left( n^{-1/d} \right)$ in \cite{KM13}.

\C{More general point processes}
\noindent
Yogeshwaran and Adler obtained similar results for Betti numbers in a much more general setting of stationary point processes \cite{YA15}.

\B{PERSISTENT HOMOLOGY}
\noindent
Bubenik and Kim studied persistent homology for i.i.d.\ random points on the circle, in the context of a larger discussion about foundations for topological statistics \cite{BK07}.

Bobrowski, Kahle, and Skraba studied maximally persistent cycles in $VR(n,r)$ and $C(n,r)$. They defined the \emph{persistence} of a cycle as $p( \sigma) = d(\sigma) / b(\sigma)$, and found a law of the iterated logarithm for maximal persistence.

\theorem{\cmrx\cite{BKS15}}
{Fix $d \ge 2$, and $1 \le i \le d-1$. Choose $n$ points i.i.d.\ uniformly randomly in the cube $[0,1]^d$.
With high probability, the maximally persistent cycle has persistence
$$\max_{\sigma} p ( \sigma) = \Theta \left( \frac{\log n}{\log \log n}  \right)^{1/i}.$$}

\conj{A law of large numbers for persistent homology.}
{$$\max_{\sigma} p ( \sigma) / \left( \frac{\log n}{\log \log n}  \right)^{1/i} \to C,$$
for some constant $C=C_{d,i}$.}

\Bnn{OTHER RESOURCES}
\noindent
For an earlier survey of Erd\H{o}s--R\'enyi based models with a focus on the cohomology-vanishing phase transition, see also \cite{Kahle14survey}. For a more comprehensive overview of random geometric complexes, see \cite{BK14}.

Several other models of random topological space have been studied. Ollivier's survey  \cite{Ollivier05} provides a comprehensive introduction to random groups, especially to Gromov's density random groups and the triangular model. Dunfield and Thurston introduced a new model of random 3-manifold \cite{DT06} which has been well studied since then.

\Bnn{RELATED CHAPTERS}

\noindent Chapter 22: Topological methods

\noindent Chapter 24: Embedding and geometric realization

\noindent Chapter 26: Persistent homology

\noindent Chapter 27: High-dimensional topological data analysis

\vspace{-1pc}

\Refh

\small
\bibliography{Krefs}
\end{document}